# Impact of Pre-Assessment and Post-Assessment in an Introductory Real Analysis Course

Chamila Malagoda Gamage[1]

Department of Mathematics, University of Florida, USA

**Abstract**

This study explores how pre- and post-assessments shape learning outcomes in an Introductory Real Analysis course. Pre-assessments act as learning roadmaps, highlighting prior knowledge and guiding student focus, while post-assessments measure growth and conceptual mastery. By analyzing student performance and feedback, we assess their impact on engagement, self-efficacy, and deeper mathematical understanding. The findings offer valuable insights for enhancing instructional strategies and fostering a more effective, student-centered learning experience in advanced mathematics.

**Keywords:** Pre-Assessment, Formative Assessment, College Mathematics, Mathematics Education, Real Analysis

## 1 Introduction

Assessment plays a pivotal role in shaping student learning experiences, particularly in mathematically rigorous courses such as Introductory Real Analysis. Traditional assessment methods often emphasize summative evaluations, such as final exams, which provide limited opportunities for students to reflect on their learning progress. In contrast, formative assessment techniques, including pre-assessments and post-assessments, have gained attention for their ability to enhance student engagement, guide instruction, and improve overall learning outcomes. Pre-assessments serve as diagnostic tools that help identify students' prior knowledge and potential misconceptions before new material is introduced. By providing an early snapshot of student understanding, instructors can tailor their teaching strategies to address specific learning gaps and optimize classroom instruction. Additionally, pre-assessments act as cognitive roadmaps, allowing students to anticipate key concepts and structure their learning effectively. Post-assessments, on the other hand, measure student progress and conceptual mastery, providing valuable feedback on the effectiveness of instructional methods and student learning trajectories.

[1] cgamage@ufl.edu

Research has demonstrated that incorporating pre- and post-assessments into undergraduate courses can enhance student motivation, self-regulation, and long-term retention of material (Rämö et al., 2023; Balan, 2012). Studies in other disciplines, such as business and finance, have shown that students who engage with pre-assessments tend to perform better on subsequent evaluations, as these tools help them focus on essential concepts and track their learning progress (Berry, 2008). However, while self-assessment practices have been explored in undergraduate mathematics education, the use of structured pre- and post-assessments remains relatively underutilized. Given the abstract and proof-intensive nature of Real Analysis, incorporating these formative assessment tools may provide a structured approach to support students' conceptual development and problem-solving skills.

This study aims to examine the impact of pre- and post-assessments on student learning outcomes in an undergraduate-level Introductory Real Analysis course. By analyzing student performance and feedback, we seek to determine whether these assessments enhance student engagement, boost self-efficacy, and foster a deeper understanding of mathematical concepts. The findings are expected to provide valuable guidance for refining instructional strategies and effectively incorporating formative assessments into the mathematics curriculum at the undergraduate level.

## 1.1 Related work in literature

In Suurtamm et al. (2016), the authors explore key issues in both large-scale and classroom-based assessments in mathematics education. They examine the purposes, traditions, and principles behind mathematics assessment, as well as the design of assessment tasks. The book addresses the challenges faced by educators, including conceptual, cultural, and political dilemmas, resistance to change, and the variability in adopting innovative assessment practices. To help teachers navigate these challenges and enhance their assessment strategies, the authors propose solutions such as sustained professional development, an emphasis on effective questioning techniques, and the integration of collaborative activities.

In Berry (2008), the author specifically examines the use of pre-tests as a non-graded assessment tool to assess students' pre-existing knowledge of a subject. The study was conducted in an undergraduate introductory corporate finance course, where students completed pre-tests before each of the eight course sections. The results indicated that pre-test scores improved as the course progressed, likely due to students engaging in the material in advance. The author also highlights several caveats, such as the potential for low pre-test scores to demoralize students, the risk of students not taking the pre-tests seriously, and the significant time commitment required from instructors. Despite these challenges, the study suggests that pre-tests can be an effective tool for enhancing student performance, as evidenced by higher final exam scores in the pre-test group, indicating that this approach helps students better master the material.

The author in Sanders (2019), provides a comprehensive overview of pre- and post-tests, explaining what they are, why they should be implemented, and their key characteristics. The article emphasizes how agencies can leverage data from these assessments to improve educational practices, highlighting the benefits of pre- and post-testing for students, teachers, and administrators. Regular data collection through these assessments plays a crucial role in monitoring the progress of youth, particularly those at risk for academic deficits, ensuring that appropriate support and programming are provided. The article also covers important considerations for selecting tests, such as ensuring their validity, choosing the right test types, and offers examples of effective pre- and post-testing practices.

## 2 Context

The context of this study is an introductory real analysis course at a public research university in the United States. The class consists of 28 students, primarily juniors and seniors, majoring in mathematics, statistics, computer science, astrophysics, physics, and economics. This course is a rigorous, proof-based mathematics class designed for students who have already completed the calculus sequence. It revisits familiar concepts from calculus such as the real number system, functions of one variable, limits, continuity, and differentiability but explores them with greater precision and mathematical rigor.

A key objective of the course is to develop students' ability to construct and communicate mathematical proofs, a skill that is often new and challenging for those accustomed to computational approaches in calculus. The course emphasizes logical reasoning, theorem formulation, and formal proof-writing techniques, requiring students to transition from procedural problem solving to abstract mathematical thinking. Instructional methods include a combination of traditional lectures and problem discussions. Students are frequently assigned proof-writing exercises and are encouraged to collaborate on challenging problems to deepen their understanding.

Common challenges in the course include difficulties in grasping formal definitions, constructing rigorous arguments, and bridging the gap between intuitive calculus concepts and their formal mathematical foundations. Many students struggle with the level of abstraction required and need time to develop confidence in their proof-writing abilities. Given these challenges, assessing students' progress in both conceptual understanding and proof-writing skills is crucial. This study investigates the effectiveness of pre- and post-assessments in measuring students' learning gains and identifying areas where additional instructional support may be needed.

## 2.1 Method

This study focused on the topic of **limits**, one of the most fundamental and challenging concepts in calculus. The activity was designed to assess students' prior knowledge, intuition, and learning gains through a structured pre- and post-assessment approach.

At the beginning of the unit, students completed a pre-assessment quiz (see Appendix A), a 15-minute diagnostic test designed to evaluate their background knowledge of limits and gauge their intuitive understanding of the topic. Immediately after completing the quiz, students self-graded their responses while the instructor explained the solutions, facilitating a discussion on key concepts. Following this, students were asked to respond to a reflection question: **"What is one mathematical concept related to limits that you would like to improve on, and what is something you would like to learn more about?"** This question encouraged students to self-assess their understanding and set learning goals for the unit.

Over the next five class sessions (each 50 minutes long), students engaged in proof-based learning of limits. The instruction covered fundamental topics such as the formal definition of a limit, one-sided limits, and limit properties. Unlike computational approaches commonly used in introductory calculus, the emphasis in this course was on proving limits rigorously rather than merely calculating them. The lectures incorporated interactive discussions and exercises in different difficulty levels to strengthen students' proof-writing skills and deepen their conceptual understanding.

At the end of the unit, students completed a post-assessment quiz (see Appendix B), designed to closely mirror the pre-assessment in both structure and content. Both quizzes contained different types of questions, including true/false, graphical interpretation, computational problems, and short proof-writing exercises. This diverse question format ensured a comprehensive evaluation of students' understanding, spanning conceptual reasoning, problem-solving skills, and formal mathematical justification. This structure allowed for a direct comparison to determine learning gains. Finally, students completed a feedback form reflecting on their experience with the pre/post-assessment process, providing insights into their learning progress and the effectiveness of the activity.

## 3 Results

Out of 28 students, 18 students took the pre-assessment quiz, and 20 students took the post-assessment quiz. To analyze the results, we computed basic statistical measures such as the mean and median. The mean pre-assessment score was **6.39**, while the mean post-assessment score increased to **11.11**, indicating a **73.91%** improvement. Similarly, the median score increased from **6.5** to **11.0**. Students' scores improved significantly after instruction. The median and mean increased notably, and the standard deviation remained similar, indicating consistent improvement. Additionally, we conducted a **paired t-test** to determine whether the improvement was statistically significant, yielding a p-value of **0.0013**, confirming that the difference was not due to random chance. To assess the magnitude of improvement, we calculated **Cohen's d effect size**, which resulted in **2.36**, indicating a large effect and strong impact of instruction. The table below summarizes the statistics.

| Measurement | Pre-Assessment | Post-Assessment |
|---|---|---|
| Mean | 6.39 | 11.11 |
| Median | 6.5 | 11.0 |
| Standard Deviation | 2.06 | 1.95 |
| Effect Size | 2.36 | |
| Paired t-test | 0.0013 | |

Table 1

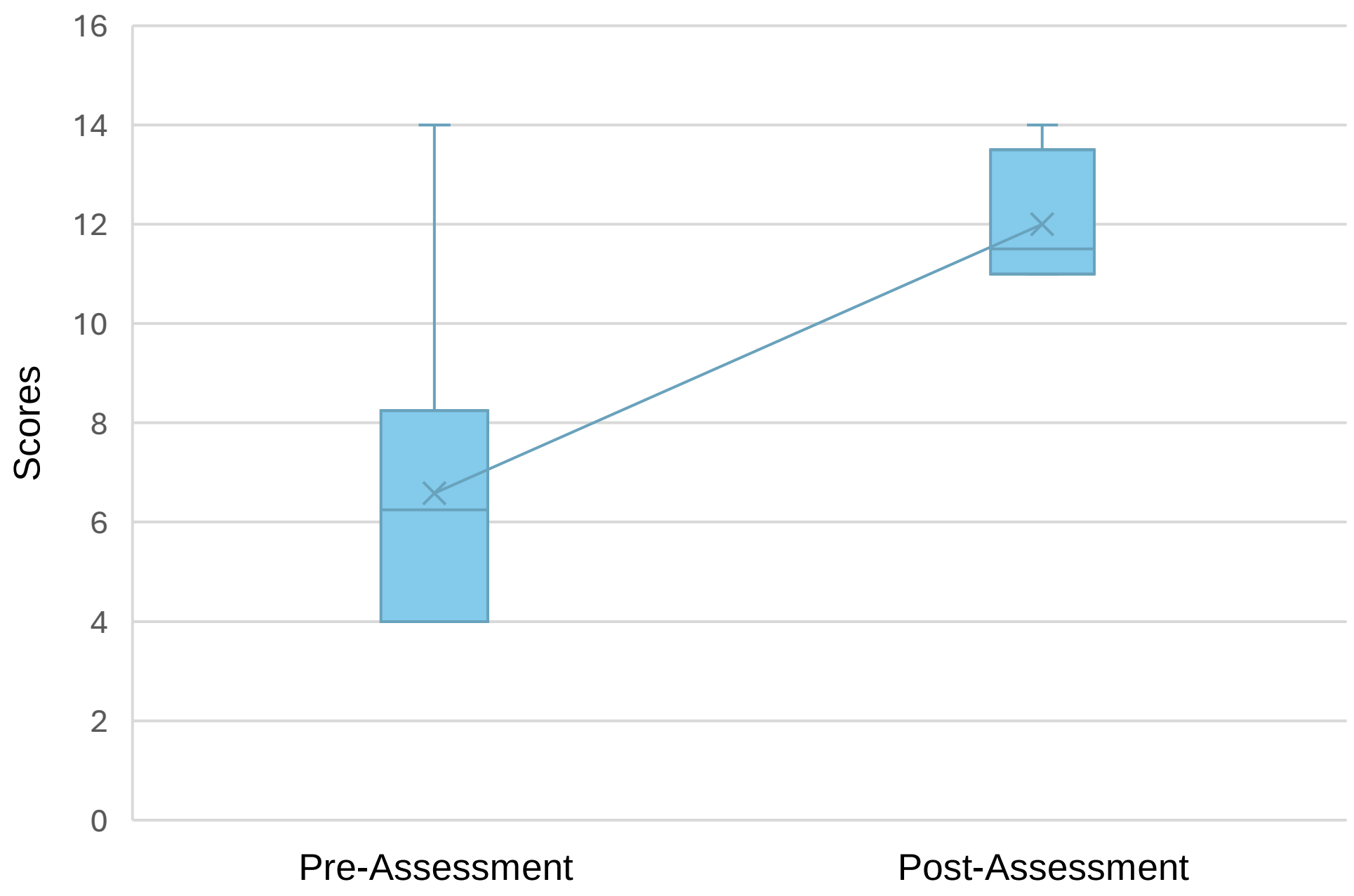

Figure 1: Pre vs. Post Assessment Scores

## 3.1 Discussion

### 3.1.1 Results from the Pre-Assessment Quiz and the Reflection Question

The pre-assessment quiz provided valuable insights into students' prior knowledge and helped shape the approach to the course. Here are some key findings:

**(I) Self-Assessment and Reflection:**

One important aspect of the pre-assessment was the requirement for students to self-grade their responses. During the subsequent solution discussion, students went beyond simply assigning points. They took the opportunity to reflect on what they didn't know, and crucially, many were able to articulate their understanding in their own words. This self-reflection is an essential aspect of fostering awareness, as students actively engaged with their gaps in knowledge and were able to communicate those gaps in a more personalized manner. This also provided the instructor with a clearer picture of where each student stood, which would be beneficial for tailoring subsequent lessons and support.

**(II) Background Knowledge Gaps:**

Another key discovery emerged from the reflection questions, which revealed that many students had forgotten important background knowledge from their previous calculus courses. Specifically, techniques related to limits, such as L'Hôpital's Rule and the Squeeze Theorem, were often unclear to students. This feedback was instrumental for the instructor, as it highlighted areas that required immediate review. In response, I shared a review workbook focusing on fundamental limit concepts, providing students with an opportunity to refresh and strengthen their foundational knowledge before delving deeper into real analysis topics.

**(III) Interest in Conceptual Understanding:**

A noteworthy observation from the pre-assessment was that many students expressed a genuine interest in understanding why limit properties work, rather than simply memorizing procedures. Students were eager to understand the actual definitions behind the theorems and the logic behind the proofs. This is encouraging, as it suggests that the class values conceptual depth and seeks to engage with the material in a meaningful way.

### 3.1.2 Results from the Post-Assessment Quiz

The post-assessment quiz demonstrated significant progress among students, with everyone showing improvement in their grades. This highlights the positive impact of the instruction and review activities throughout the course. Below are the key findings:

**(I) Improvement in Knowledge:**

One of the most encouraging results was that all students' grades increased from the preassessment to the post-assessment, indicating that their understanding of the material had deepened. Specifically, a question regarding the uniqueness of a limit; something not covered in previous calculus classes but introduced in real analysis, served as a key indicator of students' progress. In the pre-assessment, many students struggled with this concept, but by the time of the post-assessment, those who had initially answered incorrectly had now correctly answered the question, demonstrating a clearer grasp of the material.

**(II) Feedback and Reflections:**

The feedback from students, gathered through a reflection form that asked, "What was the most useful part of the post-assessment quiz for you?", was overwhelmingly positive and provided valuable insights into their learning experiences. Many students highlighted the importance of the definitions of limits and the theory behind them. Responses such as *"The definitions, because it made me realize I need to memorize them"* and *"Rewriting the definition of a limit and its existence"* underscored the significance of having a solid grasp of foundational concepts in real analysis. This suggests that students recognized the value of mastering definitions early on to build a deeper understanding of more complex topics.

**(III) Value of Reviewing and Explanation:**

Students also appreciated the opportunity to review the quiz answers, and the explanation provided during the class discussions. Several responses indicated that revisiting the questions and understanding the rationale behind the correct answers reinforced their learning. One student remarked, *"When the instructor went over the quiz and explained each* *correct answer, I understood them even more since I basically did the questions twice (once by myself and once with the instructor).”* This feedback emphasizes the benefits of engaging in active learning through reflection and instructor-guided clarification.

**(IV) Sense of Progress and Confidence:**

Many students found the post-assessment to be a good measure of their progress, with some stating that it made them feel like they were *“on track.”* This sense of progress is important for student motivation and confidence. One response mentioned that the quiz was a *“good measure of my progress,”* while another noted the usefulness of the cumulative review, which helped them *“notice what I need to study.”* The quiz allowed students to see how their understanding of key concepts, particularly limits, had improved over time.

### 3.1.3 Feedback on Overall Activity

The feedback from students regarding the overall assessment activity highlighted its positive impact on helping them achieve the learning outcomes for the unit. Many students found the activity valuable in several keyways, which contributed to their understanding and preparation for the course. The following pie charts offer a visual overview of student responses for a few questions in the feedback form:

*Q1: How do you think you performed on the post-assessment quiz compared to the pre-assessment?*

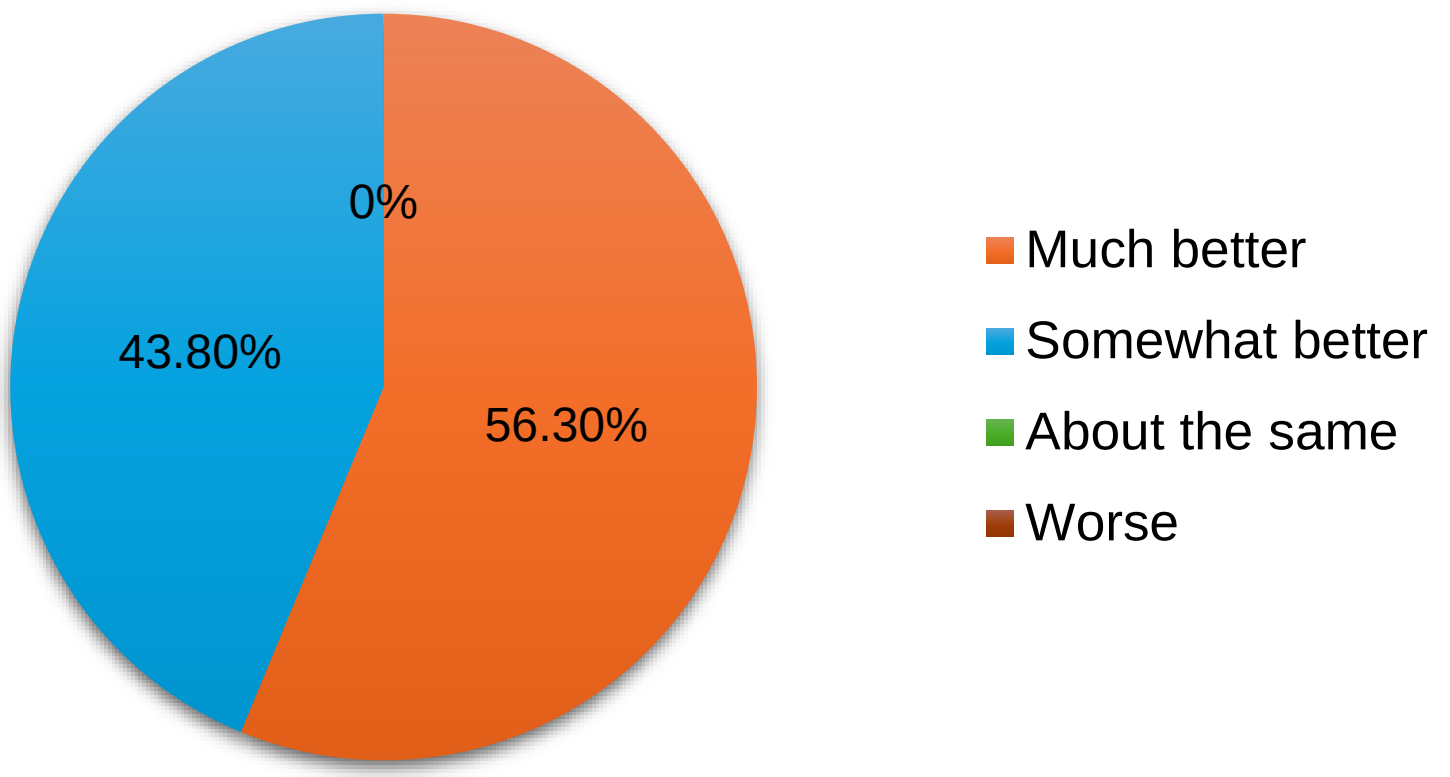

Figure 2: Performance on the post-assessment

*Q2: Did the post-assessment quiz help reinforce your understanding of limits?*

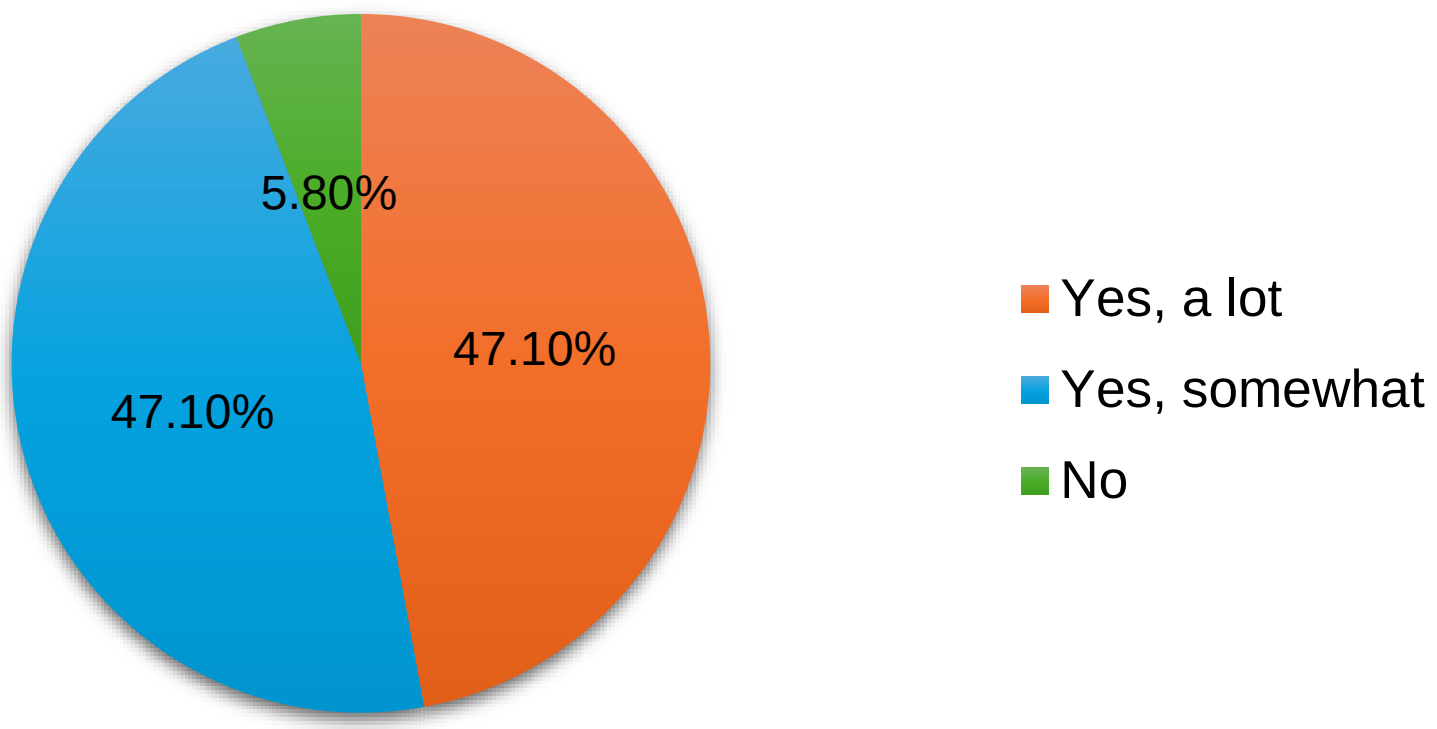

Figure 3: Effect on understanding

*Q3: Did the pre-assessment and post-assessment quizzes together help you learn better?*

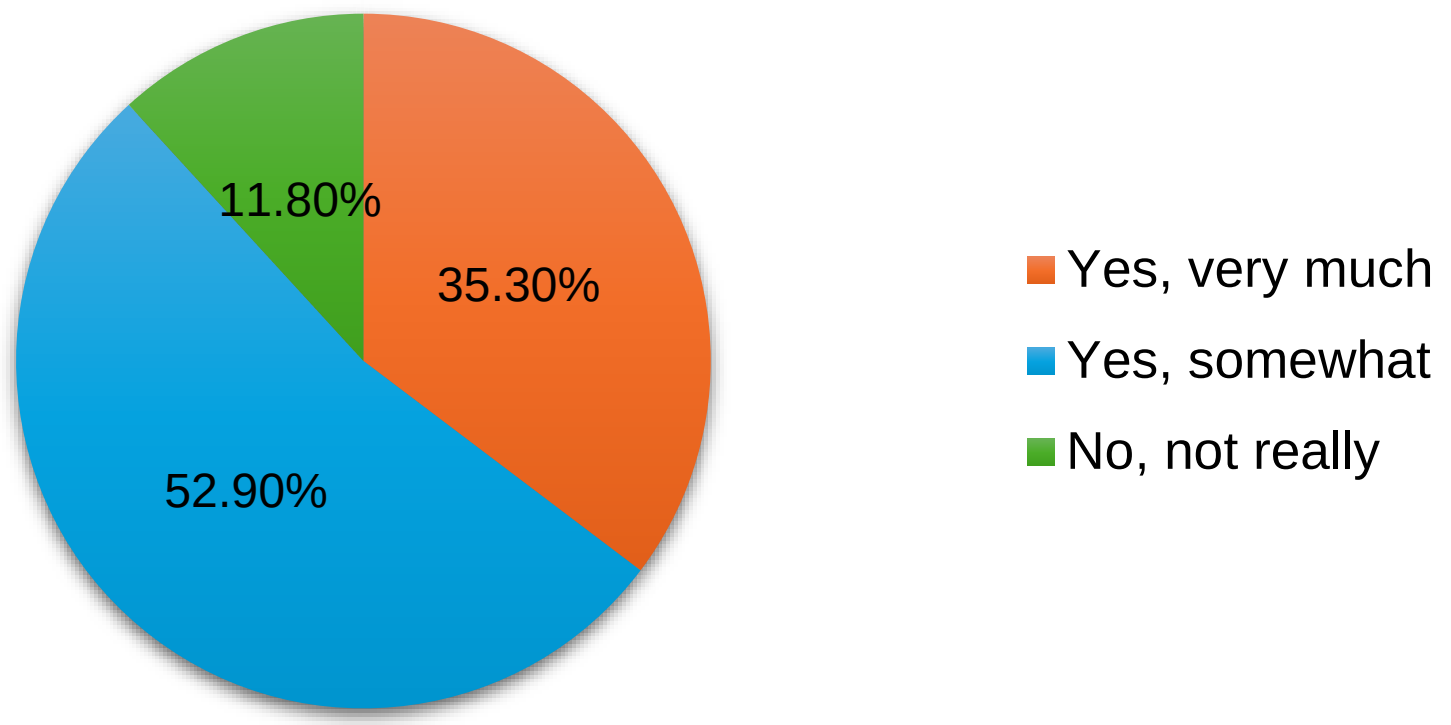

Figure 4: Feedback on the overall activity

Here is the summary of student feedback:

**(I) Preparation and Focus:**

One common theme in the feedback was that the assessment activity helped students get a clearer understanding of the course material and allowed them to focus their attention on specific areas. Many students expressed that the activity gave them a preview of what they would be learning, making it easier for them to grasp the material once it was formally introduced. As one student noted, *"It helped give me an idea of what we would be learning, which made it easier to learn the material when we did learn it."* Additionally, the activity helped them identify areas that required further attention, with comments such as *"It made me realize what I need to focus on for the exam"* and *"It helped me target specific challenge areas and focus on the necessary subjects."*

**(II) Reinforcement of Knowledge:**

Students also noted that the activity played a key role in reinforcing their knowledge. Several responses highlighted how the review and reflection process strengthened their understanding of key concepts, particularly the definitions of limits. One student stated, *"It reinforced my knowledge many times,"* while another commented, *"Reinforced limit* *definition and existence statement."* This indicates that repetition and active engagement with the material helped students solidify their grasp of the subject matter.

**(III) Self-Assessment and Reflection:**

The assessment activity also provided students with an opportunity to evaluate their own understanding of the material. Many students found this reflective process valuable, as it allowed them to gauge where they stood in terms of comprehension and identify areas for improvement. For example, one student shared, *"it helped me see where I stood on the material and gave me an idea of what to study more"*. Another mentioned, *"I could compare how I wrote my definitions/answers compared to the instructor"*, which indicates that the comparison between their own work and the instructor's expectations was helpful for aligning their understanding.

**(IV) Exam Preparation:**

The activity was seen as an effective tool for preparing students for upcoming assessments. Multiple responses indicated that it helped students prepare for the exam by revisiting foundational concepts. One student said, *"It helped me go through all the basic things I need to know and prepare for the exam,"* while another stated, *"It tested my knowledge and reinforced limit definitions."* These comments suggest that the activity not only tested their current understanding but also served as a useful review mechanism for exam readiness.

**(V) Building Confidence and Understanding Expectations:**

Finally, several students expressed that the activity helped them feel more confident in their understanding of the course and its requirements. For instance, one response stated, *"I feel like I have a good understanding of what is expected of me."* This confidence likely stems from the clarity provided by the assessment activity, which allowed students to see their progress and recognize the areas where they had gained proficiency.

## 4 Caveats and Suggestions for Improvement

While the assessment activity was generally well-received, there are several areas for potential improvement. There were suggestions regarding the timing and communication of the assessment activity. Some felt it would be helpful to give a heads-up on when the quizzes would take place, as this might encourage more participation. However, informing students in advance might skew the data, especially since the activity doesn’t count toward the course grade, potentially causing some students not to take it as seriously. To address this, it is important to emphasize the value and importance of the activity to ensure full engagement. As discussed in Berry (2008), the time factor is always a challenge, and integrating the quiz into the lesson itself could help. Using digital tools such as interactive quiz platforms or automated survey tools could save time and streamline the process, allowing for quicker feedback and more interactive participation. Additionally, the activity didn’t equally accommodate absent students, and this is something that needs to be addressed. Offering the quiz online would ensure all students have equal access, regardless of attendance, and help create a more inclusive experience for everyone.

### 4.1 Connection to Formative Assessment

This research highlights the important role of formative assessment in fostering student learning and guiding instructional practices. Pre-assessments serve as diagnostic tools, helping instructors identify students’ prior knowledge and potential misconceptions before introducing new content, which ensures that lessons are tailored to address the students’ specific needs. Reflection questions, an essential component of formative assessment, prompt students to self-assess their understanding, encouraging metacognitive practices and promoting active engagement in their learning process. Additionally, post-assessments provide valuable insights into students’ conceptual development, allowing instructors to assess the effectiveness of their teaching strategies. By regularly incorporating formative assessments, instructors create an environment that not only monitors learning progress but also enhances student engagement, encouraging active participation and a deeper understanding of the material.

## 4.2 Impact

The findings of this research have a significant impact on the role of educators, particularly in enhancing assessment strategies. By incorporating pre- and post-assessments into teaching, student progress can be more effectively tracked, and learning gaps can be addressed. Analyzing assessment results allows for the tailoring of instructional approaches to meet specific student needs, ensuring that teaching methods are data-driven and student-centered. Additionally, the use of self-assessment and reflection activities fosters a more active and self-regulated learning environment, encouraging students to take ownership of their learning. Understanding common misconceptions and struggles enables targeted support, particularly for students transitioning to proof-based mathematics. The research also reinforces the importance of emphasizing conceptual understanding and logical reasoning, influencing curriculum design and teaching approaches. Furthermore, sharing these findings with colleagues contributes to professional development and collaborative discussions, promoting a community of best practices for improving mathematics education.

## A Pre-Assessment Quiz

MAA 4102 Pre-Assessment Quiz - Limits Feb 10, 2025

**Name:** ____________________

1. (3 points) Fill in the blanks.

   Suppose that $f(x)$ is defined when $x$ is near the number $a$. We say that **the limit of** $f(x)$ **as** $x$ **approaches** $a$ **is** $L$, if ____________________

   ____________________.

2. (2 points) Circle the graphs that satisfy $\lim_{x \to a} f(x) = L$.

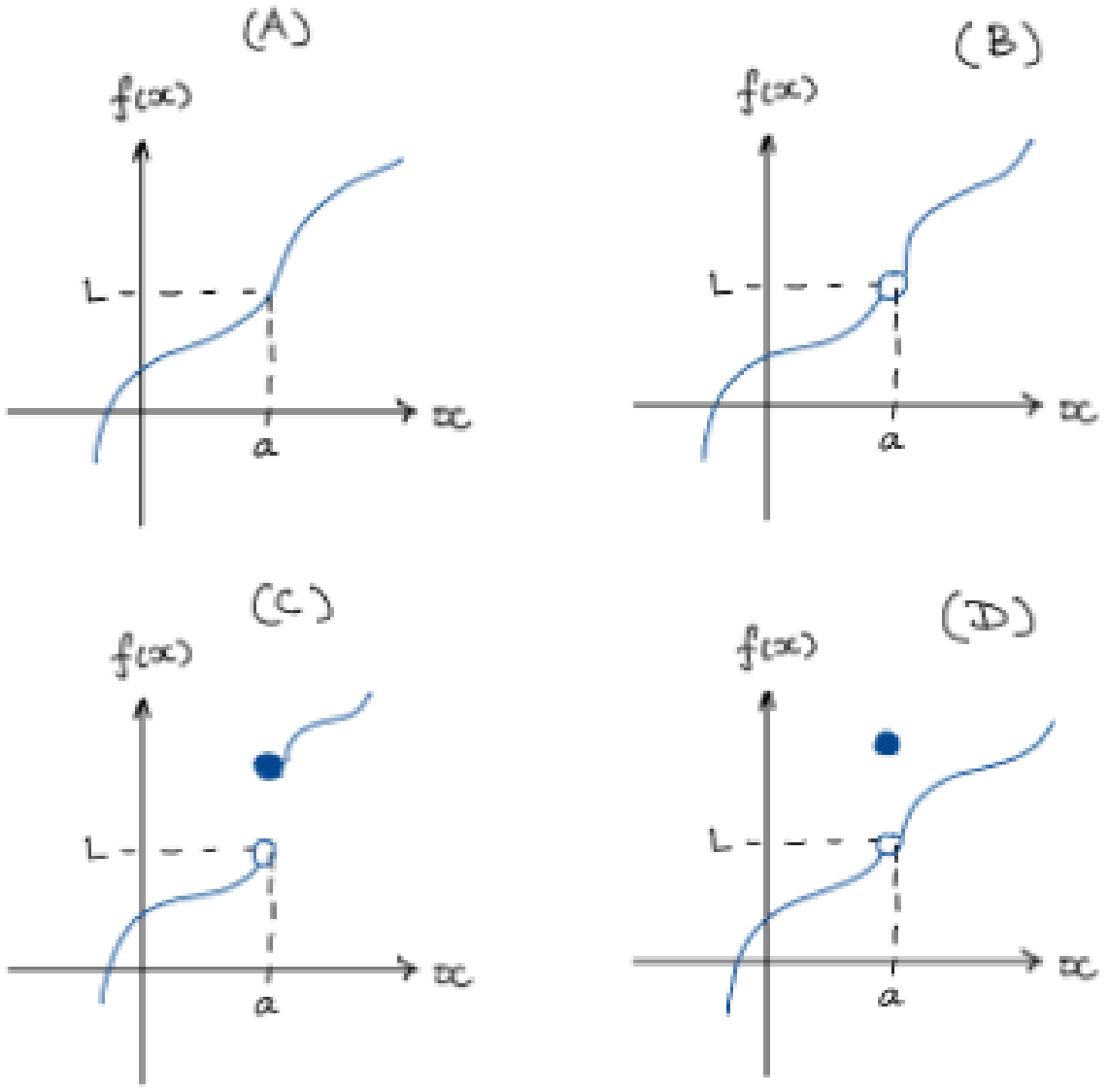

3. (1 point) Consider the function $f(x) = \begin{cases} 0 & , \quad x < 0 \\ 1 & , \quad x \geq 0 \end{cases}$.

   What is $\lim_{x \to 0} f(x)$?

4. (8 points) Evaluate the following limits. (*Show all your work!*)

   (a) $\lim\limits_{x \to 0} \dfrac{2x^2 - x}{3 + x}$

   (b) $\lim\limits_{x \to 5} \dfrac{x + 1}{x - 5}$

   (c) $\lim\limits_{h \to 0} \dfrac{\sqrt{9 + h} - 3}{h}$

   (d) $\lim\limits_{x \to 0} \left[ x^2 \cos\left(\dfrac{2}{x}\right) \right]$

5. (1 point) Let $f$ be a function such that $\lim\limits_{x \to a} f(x)$ exists. Then $\lim\limits_{x \to a} f(x)$ can take multiple values. True or False?

**Reflection:** After reviewing your quiz, what is one math concept related to 'Limits' that you would like to improve on, and what is something you would like to learn more about?

## B Post-Assessment Quiz

MAA 4102 Post-Assessment Quiz - Limits Feb 21, 2025

**Name:** ______________________________

1. (2 points) Fill in the blanks.

   Suppose that $f(x)$ is defined when $x$ is near the number $a$. We say that **the limit of** $f(x)$ **as** $x$ **approaches** $a$ **is** $L$, if ________________________________________

   ________________________________________________________________.

2. (2 points) Fill in the blanks.

   Given a function $f$, $\lim_{x \to a} f(x)$ exists iff ________________________________

   ________________________________________________________________.

3. (2 points) Circle the graphs that satisfy $\lim_{x \to a} f(x) = L$.

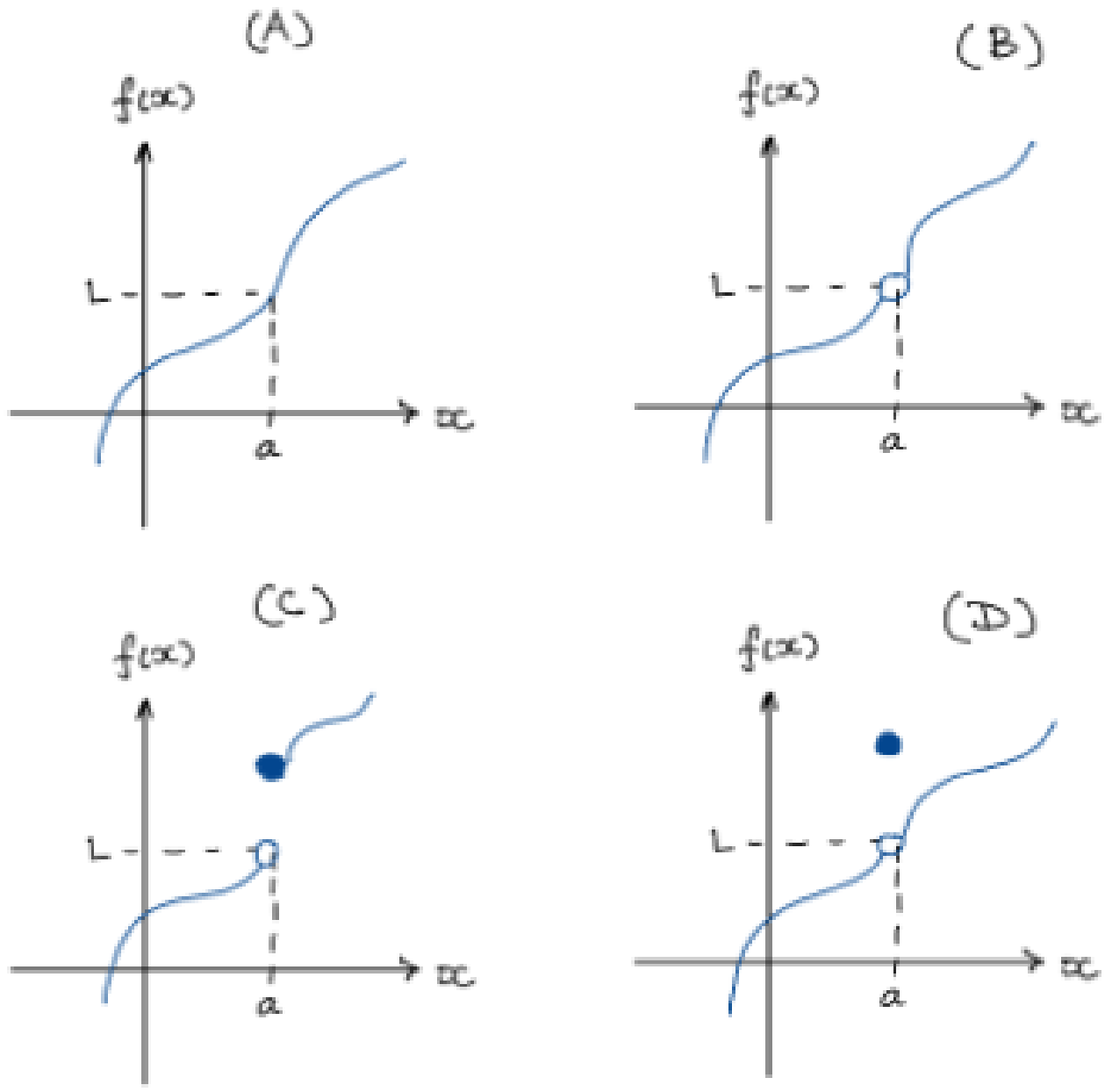

4. (1 point) Let $f$ be a function such that $\lim_{x \to a} f(x)$ exists. Then $\lim_{x \to a} f(x)$ can take multiple values. True or False?

5. (3 point) Consider the function $f(x) = \begin{cases} 2x+1 & , \quad x < 0 \\ 3 & , \quad x = 0 \\ x^2 & , \quad x > 0 \end{cases}$.

   Determine:

   (a) $\lim_{x \to 0^-} f(x)$

   (b) $\lim_{x \to 0^+} f(x)$

   (c) $\lim_{x \to 0} f(x)$

6. (1 point) If $\lim_{x \to 2} f(x) = 3$ and $\lim_{x \to 2} g(x) = -2$, then find $\lim_{x \to 2} [5f(x) - xg(x)]$.

7. (4 points) Let $\epsilon > 0$ be arbitrary. Find a $\delta > 0$ such that for all $x$ satisfying $0 < |x - 2| < \delta$, we have $|x^2 - 4| < \epsilon$.